\documentclass[10 pt,leqno]{amsart}

\usepackage {amsfonts}
\usepackage{amsthm}
\usepackage{amssymb}
\usepackage{latexsym}
\usepackage{amsmath}
\usepackage{mathrsfs}

\theoremstyle{plain}
\newtheorem{tw}{Theorem}

\newtheorem {lem} [tw]{Lemma}

\theoremstyle{definition}
\newtheorem {deft}[tw] {Definition}

\newcommand{\eqn}{\begin{equation}}
\newcommand{\eqne}{\end{equation}}

\newcommand{\bc}{\Bbb C}
\newcommand{\bn}{\Bbb N}
\newcommand{\br}{\Bbb R}

\newcommand{\cb}{\textrm{cb}}

\newcommand{\Com}{\Delta}
\newcommand{\Cou}{\epsilon}

\newcommand{\alg} {\mathsf{A}}
\newcommand{\blg} {\mathsf{B}}

\newcommand{\Tgroup}{\mathsf{T}}
\newcommand{\Ngroup}{\mathsf{N}}

\newcommand{\hil}{\mathsf{h}}

\newcommand{\kil}{\mathsf{k}}
\newcommand{\kilhat}{\hat{\mathsf{k}}}

\newcommand{\la}{\langle}
\newcommand{\ra}{\rangle}

\newcommand{\embedhn}{\iota^{(h)}_n}
\newcommand{\embedhfrac}{\iota^{(h)}_{[\frac{t}{h}]}}
\newcommand{\hfrac}{[\frac{t}{h}]}

\newcommand{\ol}{\overline}

\newcommand{\ot}{\otimes}
\newcommand{\wot}{\ol{\otimes}}
\newcommand{\ida}{\textrm{id}_{\alg}}
\newcommand{\id}{\textrm{id}}
\newcommand{\wt}{\widetilde}

\newcommand{\ve}{\varepsilon}
\newcommand{\Fock}{\mathcal F}

\newcommand{\Proc}{\mathbb{P}}

\begin{document}


\author{Uwe Franz}
\author{Adam Skalski}
\address{D\'epartement de math\'ematiques de Besan\c{c}on,
Universit\'e de Franche-Comt\'e 16, route de Gray, 25 030
Besan\c{c}on cedex, France}
\curraddr{Graduate School of Information
Sciences, Tohoku University, Sendai 980-8579, Japan}
\urladdr{http://www-math.univ-fcomte.fr/pp\underline{ }Annu/UFRANZ/}
\email{uwe.franz at univ-fcomte.fr}
\address{Department of Mathematics, University of \L\'{o}d\'{z}, ul.
Banacha 22, 90-238 \L\'{o}d\'{z}, Poland.}
\curraddr{School of Mathematical Sciences, University of Nottingham,
University Park, Nottingham, NG7 2RD }
\email{adam.skalski at nottingham.ac.uk}
\thanks{U.F.\ was supported by a Marie Curie Outgoing International
Fellowship of the EU (Contract Q-MALL MOIF-CT-2006-022137) and a Polonium cooperation.}

\title{\bf Approximation of quantum L\'evy processes by quantum random walks}

\begin{abstract}
\noindent Every quantum L\'evy process with a bounded stochastic generator is shown to arise as a
strong limit of a family of suitably scaled quantum random walks.
\end{abstract}

\maketitle

The note is concerned with investigating convergence of random walks on
quantum groups to quantum L\'evy processes. The theory of the latter is a natural
noncommutative counterpart of the theory of classical L\'evy processes on groups
(\cite{Heyer}). It has been initiated in \cite{asw} and further extensively developed
by Sch\"urmann, Schott and the first named author  (\cite{schu}, \cite{FranzSchott}, \cite{franz}).
In the series of recent papers ([LS$_{1-2}$], \cite{thesis}) Lindsay and the second named author introduced
and investigated the corresponding notion in  the topological context of compact quantum
groups (or, more generally, operator space  coalgebras).
Recent years brought also  rapid development of the theory of random walks (discrete time
stochastic processes) on discrete
quantum groups (\cite{Izumi}, \cite{Nesh}, \cite{Col}) initiated by Biane ([Bi$_{1-3}$]).

In the context of quantum stochastic cocycles (\cite{lect} and references therein)
the approximation of continous time evolutions by random walks was first
investigated by Lindsay and Parthasarathy (\cite{LiP}). They proved that
under suitable assumptions scaled random walks converge weakly to
$^*$-homomorphic quantum stochastic cocycles. Recently certain results on the strong convergence
have been obtained in papers \cite{Sinha} and \cite{Sahu} (see also \cite{Alex} for the
thorough analysis of the case of the vacuum adapted cocycles).
Here we  apply the ideas of the latter papers
to the approximation of quantum L\'evy processes
(continuous time processes) on a compact quantum semigroup $\alg$ by quantum random walks
(discrete time processes) on $\alg$.

\subsection*{Quantum random walks on $C^*$-bialgebras}

We start with the discussion of a notion of random walks on compact quantum semigroups.
The class contains finite quantum groups, so we are in a natural way generalising the notion of quantum random
walks considered in \cite{rand}. Here, and in everything that follows, $\ot$ denotes the spatial tensor product
of operator spaces (so in particular, also $C^*$-algebras).

\begin{deft}
A unital $C^*$-algebra $\alg$ is a $C^*$-bialgebra if it is equipped with two unital *-homomorphisms
$\Com: \alg \to \alg \ot \alg$ and $\Cou:\alg \to \bc$ satisfying the coassociativity and counit conditions:
\[ (\Com \ot \ida) \Com = (\ida \ot \Com)\Com,\]
\[  (\Cou \ot \ida) \Com = (\ida \ot \Cou)\Com = \ida.\]
\end{deft}

Fix for the rest of the note a $C^*$-bialgebra $\alg$.

\begin{deft}
Let $\blg$ be a unital $C^*$-algebra. A family of unital $^*$-homomorphisms $J_n:\alg \to \blg^{\ot n}$ ($n \in \bn_0$) is called a quantum random walk on $\alg$ with values in $\blg$ if
\[ J_0 = \Cou, \;\; J_n = (J_{n-1} \ot J)  \Com, \;\;n \in \bn.\]
\end{deft}

If $(J_n)_{n \in \bn}$ is a quantum random walk in the above sense, the family $(\wt{J}_n)_{n \in \bn}$
given by $\wt{J}_n = (\ida \ot J_n) \Com : \alg \to \alg \ot \blg^{\ot n}$ is a \emph{quantum random walk in the sense of \cite{LiP}}.
For any state $\phi$ on $\blg$ the family $(\kappa_n)_{n \in \bn_0}$ of states on $\alg$ defined by
\[ \kappa_n = \phi^{\ot n} \circ J_n, \;\; n \in \bn_0,\]
is a (discrete) convolution semigroup of states on $\alg$.

\subsection*{Main result}

We need first to establish some notations. Fix a Hilbert space $\kil$ and denote by
$\kilhat$ the Hilbert space $\bc \oplus \kil$ (sometimes written as $\bc \Omega \oplus \kil$).
We use the Dirac notation, so that for example $|\kil\ra$ denotes the space of all linear maps from $\bc$ to $\kil$.
The symmetric Fock space
over $L^2(\br_+;\kil)$ is denoted by $\Fock$ and its exponential vectors
by $\ve(f)$, whenever $f \in L^2(\br_+;\kil)$.  The usual shift semigroup of endomorphisms on
$B(\Fock)$
will be written $\{\sigma_s: s \geq 0\}$
and by $\Proc(\Fock)$ is meant the space of bounded adapted operator-valued processes on $\Fock$.
By a Fock space quantum L\'evy process $l \in \Proc (\alg; \Fock)$ is understood a $^*$-homomorphic map $l:\alg \to \Proc(\Fock)$, such that
\[ l_0(a) = \Cou(a) I_{\Fock}, \;\;\; l_{s+t}(a) = (l_s \ot (\sigma_s \circ l_t)) \Com (a)\]
($a \in \alg, s, t \geq 0$). It is said to be Markov regular if its Markov convolution semigroup
of states  $\{P_t = \la \ve(0), l_t(\cdot) \ve(0)\ra: t \geq 0\}$ is norm continuous.
For more information on quantum L\'evy processes on $C^*$-bialgebras we refer to \cite{LSqscc2}.

The proof of the main theorem is based on the following lemma.

\begin{lem} \label{GNS}
Assume that $\nu:\alg \to B(\kil)$ is a unital representation, and $\wt{\xi} \in \kil$ is nonzero. Let
$\gamma:\alg \to \bc$ and $\delta:\alg \to |\kil\ra$ be given by
\[ \gamma(a) = \langle \wt{\xi}, (\nu(a) - \Cou(a)) \wt{\xi}\rangle, \]
\[ \delta (a) = |(\nu(a) - \Cou(a)) \wt{\xi} \ra\]
($a \in \alg$). Define the map $\varphi:\alg \to B(\kilhat)$ by
\begin{equation} \varphi = \begin{bmatrix} \gamma  & \delta^{\dagger} \\ \delta  & \nu(\cdot) - \Cou(\cdot) I_{\kil}\end{bmatrix} \label{stochgen} \end{equation}

Put $\lambda = \|\wt{\xi}\|^2$. For each $h \in (0,\lambda^{-1}]$ there exists a unital $^*$-homomorphism
$\beta^{(h)}:\alg \to B(\kilhat)$ such that
\begin{equation} \beta^{(h)} = \begin{bmatrix} \beta^{(h)}_1  & \beta^{(h)}_2 \\ \beta^{(h)}_3  & \beta^{(h)}_4 \end{bmatrix}, \label{bgen}\end{equation}
$\beta^{(h)}_1: \alg \to \bc$, $\beta^{(h)}_3 = (\beta^{(h)}_2)^{\dagger}: \alg \to |\kil\ra$, $\beta_4^{(h)}: \alg \to B(\kil)$ and
for some constant $M>0$
\begin{align*}
&\| \beta^{(h)}_1 - (\Cou + h\gamma)\|_{\cb} \leq M h^2, \\
&\| \beta^{(h)}_3 - \sqrt{h} \delta\|_{\cb} \leq M h^{\frac{3}{2}}, \\
&\| \beta^{(h)}_4 - \nu   \|_{\cb} \leq M h.
\end{align*}

\end{lem}

\begin{proof}
Let $\xi= \wt{\xi} \|\wt{\xi}\|^{-1}$ and denote by $\nu_{\xi}$ the functional given by
\[ \nu_{\xi} (a) = \la \xi, \nu(a) \xi \ra, \;\;\; a \in \alg.\]

To construct the required $^*$-homomorphism let for $a \in \alg$
\begin{align*}  \beta^{(h)}_1 (a) = & (\Cou + h\gamma) (a) = (1 - \lambda h) \Cou(a) + \lambda h \nu_{\xi} (a),\\
 \beta^{(h)}_3 (a) = & \big| \sqrt{\lambda h} \left( \nu(a) \xi - \sqrt{1- \lambda h}\Cou(a) \xi +
( \sqrt{1- \lambda h} -1 ) \nu_{\xi} (a) \xi \right) \big\ra,\\
\beta^{(h)}_4 (a) = & \left( \lambda h \Cou(a)  + (2- 2\sqrt{1- \lambda h} - \lambda h) \nu_{\xi} (a)\right) | \xi \ra \la \xi | + \\
&( \sqrt{1- \lambda h} -1 ) |\nu (a) \xi \ra \la \xi | +
( \sqrt{1- \lambda h} -1 ) | \xi \ra \la \nu (a^*)\xi |
+ \nu(a),\end{align*}
where  Dirac notation has been again used.
Further let $\beta^{(h)}_2 = (\beta^{(h)}_3)^{\dagger}$ and define $\beta^{(h)}$ as the matrix \eqref{bgen}.

It may be checked that $\beta^{(h)}$ satisfies all requirements of the lemma.
\end{proof}

Some remarks are in place. In fact $\beta^{(h)}$ in the proof above has been constructed via the GNS construction for the state $(1 - \lambda h) \Cou + \lambda h \nu_{\xi}$. The GNS triple may be realised by  $(\Cou \oplus \nu, \bc \Omega \oplus \kil, \Omega_h)$, where
$\Omega_h = \sqrt{1 - \lambda h} \Omega \oplus  \sqrt{\lambda h} \xi$. Defining $\wt{\beta} = \Cou \oplus \nu$,
\begin{align*} \ol{\beta}^{(h)}_1 (a) & = P_{\bc \Omega_h} \ol{\beta}^{(h)} (a) P_{\bc \Omega_h}, \\
\ol{\beta}^{(h)}_2 (a) &= P_{\bc \Omega_h} \ol{\beta}^{(h)} (a) P_{(\bc \Omega_h)^{\perp}}, \\
\ol{\beta}^{(h)}_3 (a) &= P_{(\bc \Omega_h)^{\perp}} \ol{\beta}^{(h)} (a) P_{\bc \Omega_h}, \\
\ol{\beta}^{(h)}_4 (a) &= P_{(\bc \Omega_h)^{\perp}} \ol{\beta}^{(h)} (a) P_{(\bc \Omega_h)^{\perp}} \end{align*}
$(a \in \alg)$, it remains to `rotate' the GNS space to $\kilhat$ so that the decomposition
$\bc \Omega_h \oplus (\bc \Omega_h)^{\perp}$ corresponds to $\bc \Omega \oplus \kil$.
This is achieved by applying the unitary
$U_h: \bc \Omega_h \oplus \kil \to \bc \Omega \oplus \kil$ given by
\[ U_h(\alpha \Omega_h \oplus \alpha' \Sigma_h \oplus \eta) =
     \alpha \Omega \oplus \alpha' \xi \oplus \eta,\]
where $\Sigma_h = - \sqrt{\lambda h} \Omega_h \oplus  \sqrt{1-\lambda h} \xi$ and $\eta
\in (\bc \Omega_h)^{\perp} \cap (\bc \Sigma_h)^{\perp}$. It remains to check that
the maps given by ($a \in \alg$)
\begin{align*} \beta^{(h)}_1 (a) &: = U_h \ol{\beta}^{(h)}_1 (a) U_h^*, \\
\beta^{(h)}_2 (a) &: = U_h \ol{\beta}^{(h)}_2 (a) U_h^*,\\
\beta^{(h)}_3 (a) &: = U_h \ol{\beta}^{(h)}_3 (a) U_h^*,\\
\beta^{(h)}_4 (a) &: = U_h \ol{\beta}^{(h)}_4 (a) U_h^*,
\end{align*}
reduce indeed to the ones given by  formulas in the proof above. This can be done via straightforward (though very tedious)
 calculations. Note that then the fact that $\beta^{(h)}$ is a unital $^*$-homomorphism
 follows immediately from the analogous property of $\Cou \oplus \nu$. We suggest to the reader that it is worth to analyse carefully what happens to each part of the above construction as $h$ tends to $0$.
Note also that the construction of $\beta^{(h)}$ with all the properties formulated in the lemma becomes trivial if $\wt{\xi}=0$.

We are now ready to formulate and prove the main theorem of the paper.

\begin{tw}
Let $l \in \Proc (\alg; \Fock)$ be a Markov regular Fock space quantum L\'evy process on $\alg$. There exists a family of quantum random walks $(J_n^{(h)})_{n \in \bn_0}$ on $\alg$ with values in $B(\kilhat)$ and a
family of injective embeddings $\embedhn: B(\kilhat)^{\ot n} \hookrightarrow B(\Fock)$,
given by discretisation of the Fock space, (indexed by a parameter $h \in (0,\mu]$ for some $\mu >0$)
 such that 
 for each $a \in \alg$, $t \geq 0$, $\zeta \in \Fock$,
\[ (\embedhfrac  \circ J_{[\frac{t}{h}]}^{(h)} (a) ) (\zeta) \stackrel{h \to 0^+}{\longrightarrow} l_t(a) \zeta.\]
\end{tw}
\begin{proof}
Theorem 6.2 of \cite{LSqscc2} implies that the cocycle $l$ is stochastically generated by
a map $\varphi:\alg \to B(\kilhat)$ given by the formula \eqref{stochgen} for some vector
$\wt{\xi} \in \kil$ and representation $\nu: \alg \to B(\kil)$. We may assume that the vector
$\wt{\xi}$ is nonzero; otherwise the approximation method described below still works, and there is no need to restrict the range of $h>0$ in any way (see the remark before the theorem).

Let $\lambda = \|\wt{\xi}\|^2$. Let, for each $h \in (0, \lambda^{-1}]$,
$\beta^{(h)}: \alg \to B(\kilhat)$ be a $^*$-homomorphism satisfying all the
properties described in Lemma \ref{GNS}. Define the approximating
random walk by the formulas
\[ J_0^{(h)} =\Cou, \;\;\;  J_1^{(h)}  = \beta^{(h)} , \]
\[ J_{n+1}^{(h)}(a) =  (J_{n}^{(h)} \ot J^{(h)}) \Com, n \in \bn.\]
The embeddings $\embedhn$ are given by the standard discretization procedure for the Fock space
(\cite{Sahu}, \cite{Attal}). Precisely speaking , take any
$T^{(1)}, \ldots, T^{(n)} \in B(\kilhat)$,
\[ T^{(i)} = \begin{bmatrix} T^{(i)}_1  &  T^{(i)}_2 \\ T^{(i)}_3  &  T^{(i)}_4\end{bmatrix},\]
and write $\Tgroup = T^{(1)} \ot \cdots \ot T^{(n)}$. Then
\[ \embedhn (\Tgroup) \ve(f) = \bigotimes_{i=1}^k \Ngroup^{(h)}_i (T^{(i)}) \ve(f_{[(i-1)h, ih[}) \ot
\ve(f_{[kh}),\]
where
\[ \Ngroup^{(h)}_i = \sum_{l=1}^4 N^l_{T^{(i)}_l}[(i-1)h, ih], \]
and the operators $N^l_{T^{(i)}_l}$ are discretised versions of time ($N^1$),
annihilation ($N^2$), creation  ($N^3$) and preservation ($N^4$) integral,
defined as in \cite{Sahu}.

The idea of the proof is to pull the situation back to the realm of standard Markov stochastic cocycles and apply a slightly improved version of the main theorem of \cite{Sahu}.
To this end assume that $\alg$ is faithfully and nondegenerately represented on a Hilbert space $\hil$. Define
\[ \wt{l}_t = (\ida \ot l_t) \Com: \alg \to \alg \ot B(\Fock),\]
\[\wt{J}_n^{(h)} = \left(\ida \ot (\embedhn \circ J_n^{(h)})\right) \Com: \alg \to \alg \ot B(\Fock).\]
Lemma 4.1 and Proposition 3.3 of \cite{LSqscc2} imply that $\wt{l}$ is stochastically generated by an operator
$\phi = (\ida \ot \varphi) \Com: \alg \to \alg \ot B(\kilhat)$.
It may also be shown
that $\wt{J}_n^{(h)}$ coincides with the map $p_{nh}^{(h)}:\alg \to \alg \ot B(\Fock)$ constructed as in \cite{Sahu} via the $^*$-homomorphisms
$\wt{\beta}^{(h)}= (\ida \ot \beta^{(h)}) \Com: \alg \to \alg \ot B(\kilhat)$.
It is easy to note that the conditions of the Lemma \ref{GNS} imply that
if
\[ \Phi = \begin{bmatrix} \Phi_1  &  \Phi_2 \\ \Phi_3  &  \Phi_4\end{bmatrix},\;\;\;
\wt{\beta}^{(h)} = \begin{bmatrix} \wt{\beta}^{(h)}_1  &  \wt{\beta}^{(h)}_2 \\ \wt{\beta}^{(h)}_3  &  \wt{\beta}^{(h)}_4\end{bmatrix},\]
then
\begin{align*}
&\| \wt{\beta}^{(h)}_1 - \Phi_1\|_{\cb} \leq M h^2, \\
&\| \wt{\beta}^{(h)}_2 - \sqrt{h} \Phi_2\|_{\cb} \leq M h^{\frac{3}{2}}, \\
&\| \wt{\beta}^{(h)}_3 - \sqrt{h} \Phi_3\|_{\cb} \leq M h^{\frac{3}{2}}, \\
&\| \wt{\beta}^{(h)}_4 - \Phi_{4}   \|_{\cb} \leq M h.
\end{align*}

Now one may check that this is sufficient for all the assumptions of the principal theorem of \cite{Sahu} to be satisfied,
and we deduce the following statement:
for each $a \in \alg, v \in \hil$ and $\zeta \in \Fock$
\[ \wt{J}_{[\frac{t}{h}]}^{(h)} (a) (v \ot \zeta) \stackrel{h \to 0^+}{\longrightarrow} \wt{l}_t(a) (v \ot \zeta).\]

The careful analysis of the estimates used in the proof of the theorem mentioned above
shows that in fact one can obtain a
stronger result, which is of use for what follows.
Define for each $\zeta \in \Fock$, $t \geq 0,$ $ n \in \bn$ the maps
$\wt{l}_{t, \zeta}: \alg \to B(\hil; \hil \ot \Fock)$ and  $\wt{J}_{n, \zeta}^{(h)}: \alg \to B(\hil; \hil \ot \Fock)$ by the formulas
\[ (\wt{l}_{t, \zeta}(a)) (v) = \wt{l}_{t} (a) (v \ot \zeta), \]
\[ (\wt{J}_{n, \zeta}^{(h)} (a)) (v) = (\wt{J}_{n}^{(h)} (a)) (v \ot  \zeta)\]
($a \in \alg, v \in \hil$).
It is easy to see that in our context both $\wt{l}_{t, \zeta}$ and $\wt{J}_{n, \zeta}^{(h)}$
take indeed values in the operator space $\alg \ot |\Fock \rangle$; in the general, von Neumann algebraic
framework of \cite{Sahu} they would take values in the von Neumann
module $\alg'' \wot |\Fock\rangle$. As all the estimates in \cite{Sahu} are independent
of $v\in \hil$, it may be deduced in fact that for each $\zeta \in \Fock$, $t \geq 0,$ $ a\in \alg$
\begin{equation} \wt{J}_{[\frac{t}{h}], \zeta}^{(h)} (a)
 \stackrel{h \to 0^+}{\longrightarrow} \wt{l}_{t, \zeta}(a).
\label{finalest}\end{equation}
Simple argument (\cite{LSqscc2}, \cite{thesis}) shows also that  for each $a \in \alg$,
$\zeta \in \Fock$, $t \geq 0,$ $ n \in \bn$
\begin{align*} l_{t} (a) \zeta =  (\Cou \ot \id_{|\Fock\ra}) \circ \wt{l}_{t, \zeta}(a), \\
 \left(\embedhn \circ J_{n}^{(h)} (a)\right) \zeta =  (\Cou \ot \id_{|\Fock\ra})
\circ \wt{J}_{n, \zeta}^{(h)} (a). \end{align*}
In conjunction with \eqref{finalest} we obtain ($a \in \alg, \zeta \in \Fock$)
\begin{align*} \| l_{t} (a) \zeta  - \left(\embedhfrac \circ J_{[\frac{t}{h}]}^{(h)} (a)\right) \zeta \| =
\| (\Cou \ot \id_{|\Fock\ra}) \left( \wt{l}_{t, \zeta}(a) -  \wt{J}_{\hfrac, \zeta}^{(h)} (a) \right)\| \leq&  \\
\|   \wt{l}_{t, \zeta}(a) -  \wt{J}_{\hfrac, \zeta}^{(h)} (a) \|&  \stackrel{h \to 0^+}{\longrightarrow}   0
\end{align*}
This ends the proof.
\end{proof}

The main theorem above could be obtained without appealing at all to the theory of
standard quantum stochastic cocycles, essentially by rewriting the proof of L.\,Sahu
replacing everywhere the composition by the convolution operation. This is possible
only in the context of completely bounded operators; consequently, the original proof
of \cite{Sahu} would have to be formulated solely in the language of the `column'
operators (an element of a reasoning of that type may be seen in the proof above).

Markov-regular Fock space quantum L\'evy processes may be thought of as compound Poisson processes (\cite{franz}, \cite{LSqscc2}). It is therefore easy to describe conceptually how our approximations are built: the quantum random walk constructed above, after embedding in the algebra of Fock space operators, corresponds to taking random jumps, governed by the generating measure of the original compound Poisson process scaled by $h$, at discrete times $h$, $2h$, etc.. It is then clear that the limit
as $h\to 0^+$ yields the original process. The case of L\'evy processes with unbounded generators is classically resolved via treating separately the part of the process responsible for `big' jumps  and the continuous/`small' jumps part (for the extensive bibliography of the subject and applications for numerical simulations of stochastic processes we refer to \cite{Gareth}); it is not clear how to apply this procedure in the noncommutative framework.

\subsection*{Acknowledgment}
The work on this paper was initiated during the visit of the second author to the Department of Mathematics of University of Besancon in May 2006. AS would like to express his gratitude to Uwe Franz and Ren\'e Schott for making this visit possible and to Quanhua Xu for his friendly welcome and many useful discussions. This work was completed while first author was visiting the Graduate School of Information Sciences of Tohoku University as Marie-Curie fellow. He would like to thank to Professor Nobuaki Obata and the members of the GSIS for their hospitality.

\end{document}